\documentclass{scrartcl}

\usepackage[utf8]{inputenc}
\usepackage[english]{babel}
\usepackage{amssymb, amsmath, amsfonts, amsthm}
\usepackage{enumerate}
\usepackage{colonequals}
\usepackage{empheq,fancybox}
\newcommand{\euler}{\mathrm{e}} 
\newcommand{\drm}{\mathrm{d}}
\newcommand{\dvol}{\mathrm{dvol}}
\newcommand{\RR}{\mathbb{R}}

\newcommand{\NN}{\mathbb{N}}

\newcommand{\cC}{\mathcal{C}}

\renewcommand{\epsilon}{\varepsilon}
\renewcommand{\phi}{\varphi}

\DeclareMathOperator{\diam}{\mathop{diam}}

\DeclareMathOperator{\Ric}{\mathop{Ric}}
\DeclareMathOperator{\Vol}{\mathop{Vol}}

\DeclareMathOperator{\spec}{spec}

\newcommand{\rk}{(\rho-K)_-}
\newcommand{\kt}{\kappa_T}

\newtheorem{thm}{Theorem}[section]

\newtheorem{cor}[thm]{Corollary}
\theoremstyle{definition}

\theoremstyle{remark}
\newtheorem{remark}[thm]{Remark}
\numberwithin{equation}{section}

\begin{document}
\title{Almost positive {R}icci curvature in {K}ato sense - an extension of Myers' theorem} 
\author{Christian Rose\footnote{Institut f\"ur Mathematik, Universit\"at Potsdam, Karl-Liebknecht-Stra\ss{}e 24-25, D-14476 Potsdam,
Germany, christian.rose@uni-potsdam.de}}
\date{\today}

\maketitle

\begin{abstract}
It is shown that if the Kato constant of the negative part of the Ricci curvature below a positive level is small, then the volume of the corresponding manifold can be bounded above in terms of the Kato constant and the total Ricci curvature. Together with the results from \cite{CarronRose-18} and \cite{CheegerGromovTaylor-82}, this yields a generalization of the famous Bonnet-Myers theorem. Connections to some earlier generalizations are discussed. 
\end{abstract}

\hfill
\section{Introduction}
The famous Bonnet-Myers theorem states that manifolds with uniformly positive Ricci curvature are compact and have finite fundamental group. In particular, an explicit upper bound on the diameter can be given. Finding versions of this result assuming weaker curvature assumptions is an active topic in Riemannian geometry. One possibility to relax the uniform lower bound on the Ricci curvature relies on some integral conditions on the Ricci curvature along geodesics to control the index form. Another is imposing some uniform lower Ricci curvature bound and additional other assumptions on the Laplacian of the distance function or the heat kernel to conclude compactness. For a non-complete overview, see, e.g., \cite{GongWang-01,MastroliaRimoldiVeronelli-12,Wu-17} and the references therein. In contrast to the above mentioned possiblities of generalization, of particular interest for us are \cite{Li-95} and \cite{Aubry-07}, where the size of the set where the Ricci curvature fails to be uniformly bounded below by a positive constant is controlled by integral conditions on subsets of positive measure. In \cite{Aubry-07} the author proves that the volume of a manifold is always bounded from above in terms of the $L^p$-mean of the Ricci curvature below a positive threshold once the latter is finite for some $p>n/2$. Moreover, if it is small enough, then the conclusions of the Bonnet-Myers theorem hold with an explicit diameter bound. In \cite{Li-95}, it is shown that lower bounded Ricci curvature and uniform stochastic positivity of the Ricci curvature imply finite volume. This together with asymptotically non-negative Ricci curvature yields compactness by \cite{CheegerGromovTaylor-82}.
Later, it was asked in \cite{LiElworthyRosenberg-93} whether spectrally positive Ricci curvature implies compactness, a condition which is weaker than stochastic positivity of the Ricci curvature. This cannot be true in general since this condition can be satisfied by the hyperbolic space. Additional assumptions have to be imposed in order to keep the influence of the Laplacian on the Ricci curvature small enough.

As a generalization of $L^p$-Ricci curvature assumptions for $p>n/2$, Kato-type conditions on the negative part of the Ricci curvature have been used to study geometric and analytic properties of manifolds, see \cite{Carron-16,CarronRose-18,Rose-16a,RoseStollmann-15,
RoseStollmann-18,RoseWei-20}. Note that the Kato condition on the negative part of the Ricci curvature below a positive threshold implies spectrally positive Ricci curvature. An example of geometrically important features of such Kato conditions is that along the K\"ahler--Ricci flow, it is not known whether $L^p$-norms of the Ricci curvature are controllable for $p>4$, while the Kato condition is \cite{TianZhang-16}. 

This short note is dedicated to the study of compactness theorems when the negative part of the Ricci curvature is in the Kato class. 

If the manifold is assumed to be compact, diameter bounds depending on the Kato condition have been derived in \cite{CarronRose-18}, but the compactness following from Kato-type curvature conditions was missing. We will prove some compactness results in the next section under some extra hypotheses that prevent the manifold from volume collapsing at infinity. 

The first step is to show that smallness of the Kato constant of the negative part of the Ricci curvature below a positive threshold implies finite volume. This is achieved from some perturbation theoretic arguments and an old idea of Bakry \cite{Bakry-86}. The latter was already used in \cite{Li-95} under the above mentioned much stronger hypotheses. In contrast to their result, we do not need the lower bound on the Ricci curvature and derive an explicit upper bound on the volume. On the other hand, finiteness of the Kato constant is not sufficient to conclude an upper bound on the volume similarly to \cite{Aubry-07}, since the perturbation theoretic methods we use here would not work. However, our volume bound is interesting since it depends on the Kato constant of the negative part of the Ricci curvature, which is weaker than $L^p$ for $p>n/2$, and the total Ricci curvature. However, Aubry's volume estimate does not follow from our result, but a good ultracontractive bound on the heat semigroup for small times in conjunction with finiteness of the $L^p$-mean of the negative part of the Ricci curvature implies our volume bound by estimating the Kato constant as in \cite{RoseStollmann-15}. Smallness of an $L^p$-mean of the negative part of the Ricci curvature $p>n/2$ implies the Kato condition if the manifold is compact \cite{RoseStollmann-15,CarronRose-18}. This is unknown assuming completeness only.\\
Once the finite volume property is proven, the announced compactness results follow from curvature conditions yielding a lower bound on the volume growth of all sufficiently large balls centered at some point, such as non-negative or even asymptotically non-negative Ricci curvature \cite{Yau-76,CheegerGromovTaylor-82}.
On the other hand,  non-collapsing of the volume of all balls of a fixed radius is sufficient as well, which can be satisfied by several analytic assumptions on the manifold. 
The bounds on the diameter and finiteness of the fundamental group follow a posteriori from \cite{CarronRose-18}.
\section{Results}
$M=(M^n,g)$ always denotes a complete Riemannian manifold of dimension $n\in\NN$ without boundary and with Ricci tensor $\Ric$. We let $$\rho(x):=\min \spec{\Ric_x}, \quad x\in M,$$
where $\Ric$ is interpreted as a pointwise endomorphism on the cotangent bundle $T_x^*M$ at $x\in M$. Abbreviating $x_-:=\max\{0,-x\}$ and $x_+:=\max\{0,x\}$ for $x\in\RR$, we assume throughout this article that  we always have $\rho_+\in L^1_{\mathrm{loc}}(M)$. The operator $\Delta\geq 0$ is the non-negative Laplace--Beltrami operator and $(P_t)_{t\geq 0}$ denotes the heat semigroup, i.e., $P_t:=\euler^{-t\Delta}$, $t\geq 0$.

For $K\geq 0$, we say that $M=(M^n,g)$ satisfies the Kato condition if there is some $T>0$ such that
$$\kt(\rk):=\left\Vert\int_0^TP_t\rk\drm t\right\Vert_\infty<1.$$

\begin{thm}\label{thm:vol} Let $K,T>0$ and assume that $$\kt((\rho-K)_-)< 1-\euler^{-KT}.$$ Then we have
\[
\Vol(M)\leq \frac1K\frac{1-\euler^{-KT}}{1-\euler^{-KT}-\kt(\rk)}\int_M\rho\ \dvol.
\]
\end{thm}

\begin{proof}
By assumption, $\kt:=\kt(\rk)<1$. Thus, we infer from \cite{Voigt-86} 
\[
\Vert \euler^{-t(\Delta-\rk)}\Vert_{\infty,\infty}\leq \left(\frac 1{1-\kt}\right)^{1+\frac tT}, \quad t>0.
\]
The Trotter-Kato product formula implies 
\[
\Vert \euler^{-t(\Delta+\rho)}\Vert_{\infty,\infty}\leq \euler^{-Kt} \Vert \euler^{-t(\Delta-\rk)}\Vert_{\infty,\infty}\leq \euler^{-Kt}\left(\frac 1{1-\kt}\right)^{1+\frac tT}, \quad t>0.
\]
Denote by  $\Delta^1\geq 0$ the Hodge--Laplacian acting on one--forms. Since $\rho_-$ is a Kato potential, we infer from \cite{Gueneysu-17a} the following semigroup domination principle:
\begin{align*}
\Vert \euler^{-t\Delta^1}\Vert_{\infty,\infty}\leq \Vert\euler^{-t(\Delta+\rho)}\Vert_{\infty,\infty}\leq \euler^{-Kt}\left(\frac 1{1-\kt}\right)^{1+\frac tT}, \quad t>0.
\end{align*}
By assumption we have 
$$K-\frac 1T\ln\left(\frac 1{1-\kt}\right)>0.$$
Hence
\begin{align}\label{thm11}
\int_0^\infty \Vert\euler^{-t\Delta^1}\Vert_{\infty,\infty}\drm t
&\leq \int_0^\infty \euler^{-Kt}\left(\frac 1{1-\kt}\right)^{1+\frac tT}\drm t\nonumber\\
&=\frac 1{1-\kt}\int_0^\infty\exp\left(-t\left(K-\frac 1T\ln\left(\frac 1{1-\kt}\right)\right)\right)\drm t\nonumber\\
&=\frac 1{1-\kt}\left(K-\frac 1T\ln\left(\frac 1{1-\kt}\right)\right)^{-1}.
\end{align}
We follow \cite{Bakry-86} to prove finite volume by contradiction. Assume that $\Vol(M)=\infty$ and let $f,g\in \cC^\infty_c(M)$. Then, we must have $P_tf\to 0$ for $t\to\infty$. Furthermore, 
\begin{align*}
\int_M g(f-P_tf)\,\dvol
&=-\int_M g\int_0^T g\partial_tP_tf\,\drm t\, \dvol
=\int_0^T\int_M g\Delta P_tf\,\dvol \, \drm t\nonumber \\
&=\int_0^T\int_M\nabla g\nabla P_tf\,\dvol \,\drm t
\leq \int_0^T \int_M\vert \nabla g\vert \vert \nabla P_tf\vert\, \dvol\, \drm t\nonumber\\
&\leq \int_0^T\Vert\nabla P_tf \Vert_\infty \drm t\,\Vert \nabla g\Vert_1
= \int_0^T \Vert \euler^{-t\Delta^1} d f\Vert_\infty \,\drm t\,\Vert\nabla g\Vert_1\nonumber \\
&\leq \int_0^\infty \Vert \euler^{-t\Delta^1}\Vert_{\infty,\infty}\,\drm t\,\Vert d f\Vert_\infty\Vert\nabla g\Vert_1.
\end{align*}
Thus, \eqref{thm11} yields
\begin{align}\label{thm2}
\int_M g(f-P_tf)\,\dvol\leq \frac 1{1-\kt}\left(K-\frac 1T\ln\left(\frac 1{1-\kt}\right)\right)^{-1}\Vert d f\Vert_\infty\Vert\nabla g\Vert_1.
\end{align}
Consider a sequence of cut-off functions $(f_n)_{n\in\NN}$ with $0\leq f_n\leq 1$ and $\Vert d f_n\Vert_\infty\leq 1/n$ which exists by completeness. Then, \eqref{thm2} implies for all $n\in\NN$ and $g\in\cC_c^\infty(M)$
$$\int_M g(P_tf_n-f_n)\dvol\leq C/n \Vert\nabla g\Vert_1.$$
Letting $t\to \infty$ and $n\to\infty$ implies
$$\int_M g\, \dvol\geq 0,$$
contradicting that $g\in\cC_c^\infty(M)$ was arbitrary. Hence, $M$ has finite volume.\\
We infer from \cite{StollmannVoigt-96, Gueneysu-17a} 
\begin{align}\label{proof2}
\rk\leq \frac{\kt}{1-\euler^{-K T}} (\Delta+K)
\end{align}
in quadratic form sense.  This yields in particular
\begin{align}\label{proof1}
\frac{\kt}{1-\euler^{-K T}} \Delta+\rho\geq K\left(1-\frac{\kt}{1-\euler^{-KT}}\right)
\end{align}
in quadratic form sense, and the right-hand side is strictly positive by assumption.
Since $\Vol(M)<\infty$, we have $0\in\spec(\Delta)$ and so the constant function $\mathbf{1}$ is in the quadratic form domain of $\Delta$. Since $\rk$ is a Kato potential, the form domain of $\rk$ contains the quadratic form domain of $\Delta$, and thus $\mathbf{1}$ is in the form domain of $\Delta+\rho$. Hence, \eqref{proof1} yields
\[
\int_M\rho\ \dvol\geq K\left(1-\frac{\kt}{1-\euler^{-KT}}\right)\Vol(M)
\]
and the claim follows. Note that Eq.~\eqref{proof2} the same reasoning also yields an interesting bound on the mean of $\rk$. 
\end{proof}


Let $x_0\in M$, $r_0(x)$ be the distance function emanating from $x_0$. A manifold is called asymptotically non-negatively Ricci-curved if there are $x_0$ and  $r_0>0$ such that we have
$$\rho(x)\geq -\frac n{n-1}\frac 1{r_0(x)^2},\quad x\in M, r\geq r_0.$$
Asymptotically non-negative Ricci curvature yields at least logarithmic volume growth of balls around $x_0$, \cite{CheegerGromovTaylor-82}.

\begin{cor}\label{main2}
Let  $T,K,r_0>0$, $n\geq 2$, $\epsilon_0:=\min\left\{\frac{4}{n+3},1-\euler^{-KT}\right\}$, and $\epsilon\in[0,\epsilon_0)$. Any complete Riemannian manifold $M$ of dimension $n$ satisfying
$$\kappa_T(\rk)\leq \epsilon$$ and either
\begin{enumerate}[(i)]
\item $M$ is asymptotically non-negatively Ricci curved for $r_0$ or
\item $v:=\inf_{x\in M} \Vol(B(x,r_0))>0$
\end{enumerate}
is compact and
\[
\diam M\leq C(n,1-\epsilon)\frac {\pi} {\sqrt{K}}.
\]
Here $$C(n,1-\epsilon)=\sqrt{\frac{1-\epsilon\frac{n-1}{n}}{1-\epsilon\frac{n+3}{4}}\left(1+\frac{\epsilon}{n-1}\right)}\to 1$$ as $\epsilon\to 0$. Moreover, the fundamental group of $M$ is finite.

\end{cor}
\begin{remark} 
The condition $v>0$ can be achieved as soon as one has a Sobolev embedding or a uniform heat kernel bound for a certain time \cite{Carron-96, SaloffCoste-02}. Thus, our results somehow fit together with those from \cite{GongWang-01}.
\end{remark}

\begin{proof}[Proof of Corollary~\ref{main2}] By Theorem~\ref{thm:vol} we have $\Vol(M)<\infty$. 
Assume that $M$ is not compact. For the first part, we infer from \cite[Theorem~4.9(ii)]{CheegerGromovTaylor-82} that an asymptotically non-negatively Ricci-curved non-compact manifold has at least logarithmic volume growth at the center point and hence must have infinite volume, so that we have a contradiction. For the second part, note that the assumption and finite volume implies volume doubling and hence at least linear volume growth \cite{Yau-76}, contradicting finite volume as well. The diameter bounds and the finiteness of the fundamental group follow from \cite[Corollary~2.4, Theorem~2.8]{CarronRose-18} by setting $k^2=\frac{K}{2(n-1)}$ in the notation of the cited paper.
\end{proof}
\begin{remark}
The proof works for any curvature condition giving finite volume and, assuming a non-compact $M$, the existence of $p\in M$, $r_0$, and $f\colon (r_0,\infty)\to (0,\infty)$ non-decreasing with $f(t)\to\infty$ as $t\to\infty$ such that $\Vol(B(p,r))\geq f(r)$, $r\geq r_0$. 
\end{remark}
\subsection*{Acknowledgements} I am grateful for the comments of the anonymous referee which helped to significantly improve the presentation of this article.
\bibliographystyle{alpha}
\def\cprime{$'$}

\end{document}